\documentclass[final]{siamltex}
\usepackage{latexsym}
\usepackage{epsfig,graphics}
\def\bd{\begin{description}}
\def\ed{\end{description}}

\def\beq{\begin{equation}}
\def\eeq{\end{equation}}
\def\bea{\begin{eqnarray}}
\def\eea{\end{eqnarray}}
\def\beas{\begin{eqnarray*}}
\def\eeas{\end{eqnarray*}}

\newtheorem{defin}{Definition}[section]

\title{ACCELERATION OF UNIVARIATE GLOBAL OPTIMIZATION ALGORITHMS
WORKING WITH LIPSCHITZ FUNCTIONS\\ AND LIPSCHITZ FIRST
DERIVATIVES\thanks{This research was partially supported by
  the project ``High accuracy supercomputations and solving global
  optimization problems using the information approach"   of
the Russian Federal Program ``Scientists and Educators in Russia of
Innovations",  project 14.B37.21.0878.}}

\author{DANIELA LERA\thanks{Dipartimento di Matematica e Informatica, Universit\`{a}
degli Studi di Cagliari, Cagliari, Italy, ({\tt
lera@unica.it}).}\and YAROSLAV D. SERGEYEV\thanks{Corresponding
author. He works at the following institutions:  University of
Calabria, Rende, Italy; N.I.~Lobatchevsky State University,
  Nizhni Novgorod, Russia; and   Institute of High Performance Computing
   and Networking of the National Research Council of
   Italy, ({\tt
yaro@si.deis.unical.it})}
   }

\begin{document}

\maketitle

\begin{abstract}
This paper deals  with two kinds of the one-dimensional global
optimization problems over a closed finite interval: (i) the
objective function $f(x)$ satisfies the Lipschitz condition with a
constant $L$; (ii) the first derivative of $f(x)$ satisfies the
Lipschitz condition with a constant $M$. In the paper, six
algorithms are presented for the case (i) and six algorithms for the
case (ii). In both cases,   auxiliary functions are constructed and
adaptively improved during the search. In the case (i),  piece-wise
linear functions are constructed and in the case (ii) smooth
piece-wise quadratic functions  are used.    The constants $L$ and
$M$ either are taken as values known a priori or are dynamically
estimated  during the search. A recent technique that adaptively
estimates the local Lipschitz constants over different zones of the
search region is used to accelerate the search. A new technique
called the \emph{local improvement} is introduced in order to
accelerate the search in both cases (i) and (ii). The algorithms are
described in a unique framework, their properties are studied from a
general viewpoint, and convergence conditions of the proposed
algorithms are given. Numerical experiments executed on 120 test
problems taken from the literature show quite a promising
performance of the new accelerating techniques.
\end{abstract}



\pagestyle{myheadings} \thispagestyle{plain} \markboth{DANIELA LERA
AND YAROSLAV D. SERGEYEV}{ACCELERATION OF UNIVARIATE GLOBAL
OPTIMIZATION ALGORITHMS}

\section{Introduction}
\setcounter{equation}{0}

Let us consider the one-dimensional global optimization problem
 of finding a point $x^*$ belonging to a finite interval $[a,b]$ and the
value $f^*=f(x^*)$ such that
\begin{equation} \label{p}
f^*=f(x^*)=\min \{ f(x): \ x \in [a,b] \},
\end{equation}
where either the objective function $f(x)$ or its first derivative
$f'(x)$ satisfy the Lipschitz condition, i.e., either
\begin{equation} \label{fun}
| f(x) - f(y) | \leq L | x-y |, \hspace{1cm}   x,y \in [a,b],
\end{equation}
or
\begin{equation} \label{der}
| f'(x) - f'(y) | \leq M | x-y |, \hspace{1cm}   x,y \in [a,b],
\end{equation}
with constants $0<L<\infty,\,\, 0< \ M<\infty$.

Problems of this kind are worthy of a great attention because of at
least two reasons. First, there exists a large number of real-life
applications where it is necessary to solve univariate global
optimization problems stated in various ways (see, e.g.,
\cite{Calvin,Calvin&Zilinskas,Casado:et:al.(2002),Daponte:et:al.(1995),Gergel(1992),Hansen(1979),
Kvasov&Sergeyev(2009),Locatelli,Piya(1972),Sergeyev(1995b),
Sergeyev(2006),Sergeyev:et:al(2001),Sergeyev&Grishagin(1994),Sergeyev&Kvasov(2008),
Sergeyev:et:al.(2007),Sergeyev&Markin(1995),Zilinskas,Wolfe}). This
kind of problems is often encountered in scientific and engineering
applications (see, e.g.,
\cite{Hamacher,kn:theory,kn:KABA89,Pinter(1996),Sergeyev(1995b),Sergeyev&Kvasov(2008),
Sergeyev:et:al.(2007),Strongin&Sergeyev(2000),Wolfe}), and,  in
particular, in electrical engineering optimization problems (see,
e.g.,
\cite{Daponte:et:al.(1995),Daponte:et:al.(1996),kn:filters,Molinaro&Sergeyev(2001b),
Sergeyev:et:al.(1999),Strongin&Sergeyev(2000)}). On the other hand,
it is important to study one-dimensional  methods proposed to solve
problems (\ref{p}), (\ref{fun}) and (\ref{p}), (\ref{der}) because
they can be successfully extended to the multi-dimensional case by
numerous schemes (see, for example, one-point based, diagonal,
simplicial, space-filling curves, and other popular approaches in
\cite{Floudas&Pardalos(1996),Horst&Pardalos(1995),Horst&Tuy(1996),Mladineo(1992),
Pinter(1996),Sergeyev&Kvasov(2008),Strongin&Sergeyev(2000)}).

\begin{figure}[t]
\centering
\centerline{\epsfxsize=12cm\mbox{\epsfbox{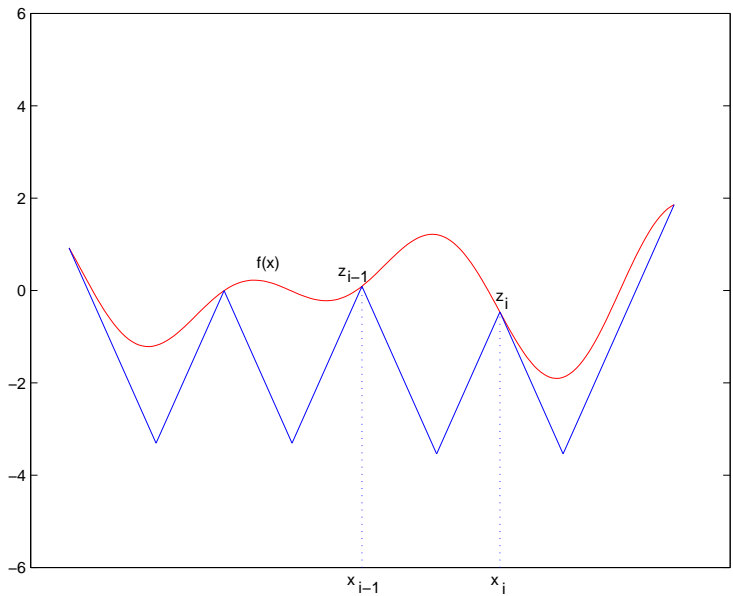}}}
\caption{\em A piece-wise linear support function constructed by
the method of Piyavskii   after five evaluations of the objective
function $f(x)$}
  \label{fig.1-2-1}
\end{figure}

In the literature, there exist several methods for solving the
problems (\ref{p}), (\ref{fun}) and (\ref{p}), (\ref{der})  (see,
for example,
\cite{Horst&Pardalos(1995),Horst&Tuy(1996),Sergeyev&Kvasov(2008),
Strongin&Sergeyev(2000),Pinter(1996)}, etc.). For solving the
problem (\ref{p}), (\ref{fun}) Piyavskii (see \cite{Piya(1972)})
has proposed a popular method that requires an a priori
overestimate of the Lipschitz constant $L$ of the function $f(x)$:
in the course of its work, the algorithm constructs piece-wise
linear support functions for $f(x)$ over every subinterval
$[x_{i-1}, x_i]$, $i=2,...,k$, where the points $x_1,...,x_k$ are
points previously produced by the algorithm (see
Fig.~\ref{fig.1-2-1}) at which the objective function $f(x)$ has
been evaluated, i.e., $z_i=f(x_i), i=2,...,k$.

In the present paper, to solve the problem (\ref{p}), (\ref{fun})
we consider  Piyavskii's method and algorithms that  dynamically
estimate the Lipschitz information for the entire region $[a,b]$
or for its subregions. This is done  since the precise information
about the value $L$ Piyavskii's method requires for its correct
work is often hard to get in practice.  Thus, we use two different
procedures to obtain an information on the constant $L$: the first
one estimates the global constant during the search (the word
``global'' means that the same estimate is used over the whole
region $[a,b]$), and the second, called ``local tuning technique''
that adaptively estimates the local Lipschitz constants in
different subintervals of the search region during the course of
the optimization process.

Then, in order to accelerate the search, we propose a new
acceleration tool, called ``local improvement'', that can be used
together with  all three ways described above   to obtain the
Lipschitz information in the framework of the Lipschitz
algorithms. The new approach forces the global optimization method
to make a local improvement of the best approximation of the
global minimum immediately after a new approximation better than
the current one is found.

The proposed local improvement technique is of a particular interest
due to the following  reasons. First, usually in the global
optimization methods the local search phases are separated from the
global ones. This means that it is necessary to introduce a rule
that  stops the global phase and starts the local one; then it stops
the local phase and starts the global one. It can happen  (see,
e.g.,
\cite{Horst&Pardalos(1995),Horst&Tuy(1996),Sergeyev&Kvasov(2008),
Strongin&Sergeyev(2000),Pinter(1996)}, etc.), that the global search
and the local one are realized  by different algorithms and the
global search is not able to use \textit{all} evaluations of $f(x)$
made during the local search losing so an important information
about the objective function that has been already obtained. The
local improvement technique introduced in this paper does not have
this defect and allows the global search to use all the information
obtained during the local phases. In addition, it can work  without
any usage of the derivatives and this is a valuable asset when one
solves the problem (\ref{p}), (\ref{fun}) because, clearly,
Lipschitz functions can be non-differentiable.

\begin{figure}[t]
\centering
\centerline{\epsfxsize=12cm\mbox{\epsfbox{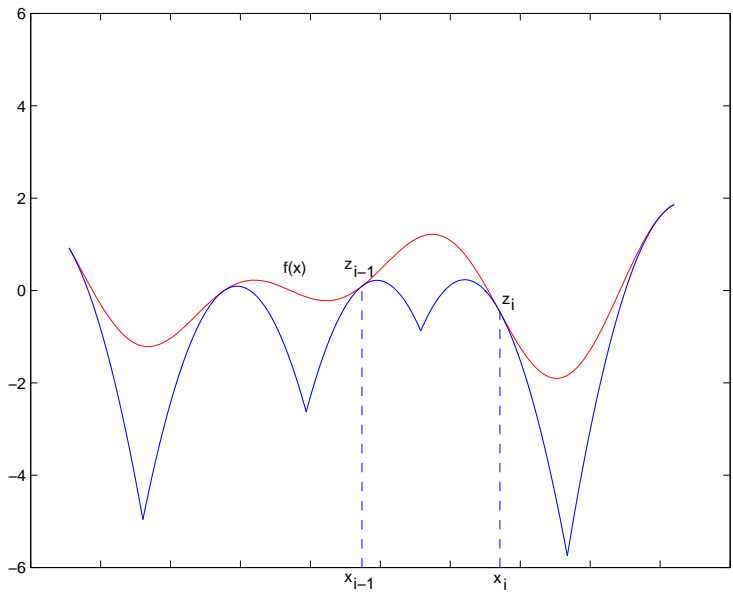}}}
\caption{\em Breiman-Cutler-Gergel piece-wise quadratic
non-differentiable support function constructed  after five
evaluations of the objective function $f(x)$}
  \label{fig.1-2-2}
\end{figure}

Let us consider now the problem (\ref{p}), (\ref{der}). For this
case, using the fact that the first derivative $f'(x)$ of the
objective function satisfies the Lipschitz condition (\ref{der}),
Breiman and Cutler (see \cite{Breiman&Cutler(1993)}) have
suggested an approach that constructs  at each iteration
piece-wise quadratic non-differentiable support functions for the
function $f(x)$ over $[a,b]$ using an a priori given overestimate
of $M$ from (\ref{der}). Gergel (see \cite{Gergel(1992)}) has
proposed independently a global optimization method that
constructs similar auxiliary functions  (see Fig.~\ref{fig.1-2-2})
and estimates $M$ dynamically during the search.

 If we suppose that the   the Lipschitz constant $M$
from (\ref{der}) is known, then (see
\cite{Breiman&Cutler(1993),Gergel(1992)}), at an iteration $k>2$,
the support functions $\Phi_i(x)$ are constructed for every
interval $[x_{i-1}, x_i]$, $i=2,...,k$, (see Fig.~\ref{fig.1-2-2})
  as follows:
\begin{equation} \label{Cutlersupp}
\Phi_i(x)=\max \{\phi_{i-1}(x), \phi_i(x)\}, \hspace{0.5cm}
x\in[x_{i-1}, x_i],
\end{equation}
where
\[
\phi_{i-1}(x)=z_{i-1}+z'_{i-1}(x-x_{i-1})-\frac{M}{2}(x-x_{i-1})^2,
\]
 \[
\phi_i(x)=z_i-z'_i(x_i-x)-\frac{M}{2}(x_i-x)^2,
\]
and $z_i=f(x_i)$, $z'_i=f'(x_i)$.

It can be noticed that in spite of the fact that $f(x)$ is smooth,
the support functions $\Phi_i(x)$ are not smooth. This defect has
been  eliminated in \cite{Sergeyev(1998)} where there have been
introduced three methods constructing smooth support functions
that are closer to the objective function than non-smooth ones.

In this paper, for solving   the problem (\ref{p}), (\ref{der}) we
describe six different algorithms where smooth support functions
are used. As it was in the case of the problem (\ref{p}),
(\ref{fun}), the local tuning and the local improvement techniques
are applied to accelerate the search.

The paper has the following  structure: in Section 2 we describe
 algorithms for solving the problem (\ref{p}), (\ref{fun}); in
Section 3 we describe methods that use smooth support functions in
order to solve the problem (\ref{p}), (\ref{der}). The convergence
conditions to the global minimizers for the introduced methods are
established in both Sections. In Section 4, numerical results are
presented and discussed. Finally, Section 5
  concludes the paper.

\section{Six methods constructing piece-wise linear auxiliary
functions for solving problems with the Lipschitz objective
function}
\setcounter{equation}{0}

In this Section, we study the problem
(\ref{p}) with the objective function $f(x)$ satisfying the
Lipschitz  condition (\ref{fun}). First, we present a general
scheme describing  in a compact form all the methods considered in
this Section and then, by specifying STEP 2 and STEP 4 of the
scheme, we introduce six different algorithms. In this Section, by
the term {\em trial} we denote the evaluation of the function
$f(x)$ at a point $x$ that is called the {\em trial point}.

 \vskip 12pt {\bf General
Scheme ($GS$) describing  algorithms working with piece-wise
linear auxiliary functions.}
 \bd
\item[STEP 0.] The first two trials are performed at the points $x^1=a$ and
$x^2=b$. The point $x^{k+1}$, $k\geq 2$, of the current (k+1)-th iteration is chosen as follows.
\item[STEP 1.] Renumber the trial points $x^1, \ x^2, \ \dots, \ x^k$ of the
previous iterations by subscripts so that
 \begin{equation} \label{order}
 a=x_1 < \dots < x_k =b.
\end{equation}
\item[STEP 2.] Compute in a certain way the values $l_i$   being
 estimates of the Lipschitz constants of $f(x)$ over the intervals
$[x_{i-1}, x_i]$, $i=2,...k$. The way  to calculate the values
$l_i$ will be specified in each concrete algorithm described
below.
 \item[STEP 3.] Calculate  for each interval $(x_{i-1},
x_i)$, $i=2,...k$,
 its characteristic
\begin{equation} \label{Mi}
 R_i= \frac{z_i + z_{i-1}}2-l_i\frac{(x_i-x_{i-1})}2,
\end{equation}
where the values $z_i=f(x_i)$, $i=1,...,k$.
\item[STEP 4.] Find an interval $(x_{t-1}, x_t)$ where the next trial will be executed.
 The way  to choose such an interval will be specified in each
concrete algorithm described below.
\item[STEP 5.] If
\begin{equation} \label{epsilon}
 |x_t - x_{t-1} | > \varepsilon,
\end{equation}
where $\varepsilon > 0$ is a given search accuracy, then execute the next trial at the point
\begin{equation} \label{point}
x^{k+1} = \frac{x_t+x_{t-1}}2+\frac{z_{t-1}-z_t}{2l_t}
\end{equation}
and go to STEP 1. Otherwise, take as an estimate of the global
minimum $f^*$ from (\ref{p}) the value
$$ f^*_k = \min \{z_i \ : \ 1\leq i \leq k \}, $$
and a   point
$$ x^*_k = \arg \min \{z_i \ : \ 1\leq i \leq k \}, $$
as an estimate of the global minimizer $x^*$, after executing
these operations STOP. \ed

Let us make   some observations with regard to the scheme $GS$
introduced above. During the course of the $(k+1)th$ iteration a
method following this scheme constructs an auxiliary piece-wise
linear function
$$ C^k(x)= \bigcup_{i=2}^{k} c_i(x) $$
where
$$ c_i(x)= \max \{z_{i-1}-l_i(x-x_{i-1}), z_i+l_i(x-x_i) \}, \hspace{10mm} x\in[x_{i-1},x_i], $$
and the characteristic $R_i$ from (\ref{Mi}) represents the
minimum of the auxiliary function $c_i(x)$ over the interval
$[x_{i-1}, x_i]$.

% and the point $x^{k+1} \in (x_{t-1}, x_t)$ is
%chosen in such a way that
%\[
%  c_t(x^{k+1}) = \min \{ C^k(x) : x \in [a,b] \} =   \min \{ R_i : 2\leq i \leq k \}.
%\]
If the constants $l_i$ are equal or larger than the Lipschitz
constant $L$ for all $i=2,...,k$, then it follows from (\ref{fun})
that the function $C^k(x)$ is a low-bounding function for $f(x)$
over the interval  $[a,b]$, i.e., for every interval
$[x_{i-1},x_i]$, $i=2,...,k$, we have
\[
f(x) \ge c_i(x), \hspace{1cm} x\in[x_{i-1},x_i],
\hspace{5mm}i=2,...,k.
\]
Moreover, if $l_i=L$, we obtain the Piyavskii support functions
(see Fig.~\ref{fig.1-2-1}).

In order to obtain from the general scheme $GS$ a concrete global
optimization algorithm, it is necessary to define   STEP~2 and
STEP~4 of the scheme. This section proposes six specific
algorithms executing this operation in different ways. In STEP~2,
we can make three different choices of computing the constant
$l_i$  that lead to  three different procedures that are called
STEP 2.1, STEP 2.2, and STEP 2.3, respectively. The first way to
define STEP 2 is the following. \vskip12pt \noindent
{\textbf{STEP 2.1.} } \\
\hspace*{1cm} Set
\begin{equation} \label{hnota}
 l_i= L, \hspace{1cm} i=2,...,k.
\end{equation}
Here the exact value of the a priori given Lipschitz constant is
used. Obviously, this rule gives us  the Piyavskii algorithm.

If the constant $L$ it is not available (this situation be can
very often encountered in practice), it is necessary to look for
an approximation of $L$ during the course of the search. Thus, as
the second way to define STEP 2 of the $GS$   we use an adaptive
estimate  of the global Lipschitz constant (see
\cite{Sergeyev&Kvasov(2008),Strongin&Sergeyev(2000)}), for each
iteration $k$. More precisely we have: \vskip12pt \noindent
{\textbf{STEP 2.2.}} \\
\hspace*{1cm} Set
\begin{equation} \label{hglob}
 l_i= r \max \{\xi, H^k\}, \hspace{1cm} i=2,...,k,
\end{equation}
\hspace*{1cm} where  $\xi >0$ is a small number that takes into account
our hypothesis that~$f(x)$  \\
\hspace*{1cm}  is not constant over the
interval $[a,b]$ and $r>1$ is a reliability parameter.~The  \\
\hspace*{1cm}  value $H^k$ is calculated as follows
\begin{equation} \label{Hmax}
H^k= \max \{H_i: i=2,...,k,\}
\end{equation}
\hspace*{1cm} with
\begin{equation} \label{Hi}
 H_i=\frac{|z_i - z_{i-1}|}{x_i - x_{i-1}},
\hspace{.5cm} i=2,...,k.
\end{equation}

In both cases, STEP 2.1 and STEP 2.2, at each iteration $k$ all
quantities $l_i$ assume the same value over the whole search
region  $[a,b]$. However, both the a priori given exact constant
$L$ and its global estimate (\ref{hglob}) can provide a poor
information about the behavior of the objective function $f(x)$
over every small subinterval $[x_{i-1},x_i] \subset [a,b]$. In
fact, when the local Lipschitz constant related to the interval
$[x_{i-1}, x_i]$ is significantly smaller than the global constant
$L$, then the methods using only this global constant or its
estimate (\ref{hglob}) can work slowly over such an interval (see
\cite{Sergeyev(1995b),Sergeyev&Kvasov(2008),Strongin&Sergeyev(2000)}).

In order to overcome this difficulty,   we consider a recent
approach (see
\cite{Sergeyev(1995b),Sergeyev&Kvasov(2008),Strongin&Sergeyev(2000)})
called the \textit{local tuning} that adaptively estimates the
values of the local Lipschitz constants related to the intervals
$[x_{i-1}, x_i],\, i=2,...,k$ (note that other techniques using
different kinds of local information in global optimization can be
found also in
\cite{Pinter(1996),Strongin&Sergeyev(2000),Torn&Zilinskas(1989)}).
The auxiliary function $C^k(x)$ is then constructed by using these
local estimates for each interval $[x_{i-1}, x_i]$, $i=2,...,k$.
This technique is described below as the rule STEP 2.3. \vskip 12pt
\noindent
{\textbf{STEP 2.3.} } \\
\hspace*{1cm} Set
\begin{equation} \label{hloc}
 l_i = r \max \{\lambda_i, \gamma_i, \xi \}
\end{equation}
\hspace*{1cm} with
\begin{equation} \label{landa}
 \lambda_i= \max\{H_{i-1}, H_i, H_{i+1}\}, \ \ i=3,...,k-1,
\end{equation}
\hspace*{1cm} where $H_i$ is from (\ref{Hi}), and when $i=2$ and
$i=k$  only $H_2,$ $H_3$,\\
\hspace*{1cm}   and $H_{k-1}, H_k$, should be considered,
respectively. The value
\begin{equation} \label{gamma}
\gamma_i = H^k \frac{(x_i-x_{i-1})}{ X^{max}},
\end{equation}
\hspace*{1cm} where $H^k$ is from (\ref{Hmax}) and
\[
 X^{max} = \max \{x_i-x_{i-1}: \ i=2,...,k \}.
\]
\hspace*{1cm} The parameter $\xi >0$  has the same sense as in
STEP 2.2. \vskip 12pt \noindent Note that in (\ref{hloc}) we
consider two different components, $\lambda_i$ and $\gamma_i$,
that take into account  respectively the local and the global
information obtained during  the previous iterations. When the
interval $[x_{i-1}, x_i]$ is large, the local information is not
reliable and the global part $\gamma_i$ has a decisive influence
on $l_i$ thanks to   (\ref{hloc})  and (\ref{gamma}). When
$[x_{i-1}, x_i]$ is small, then the local information becomes
relevant, $\gamma_i$ is small (see (\ref{gamma})), and the local
component $\lambda_i$ assumes the key role. Thus, STEP 2.3
automatically balances the global and the local information
available at the current iteration. It has been proved for a
number of global optimization algorithms that the usage of the
local tuning can accelerate the search significantly (see
\cite{Sergeyev(1995b),Sergeyev(1998),Sergeyev(2006),Sergeyev&Kvasov(2008),Sergeyev:et:al.(2007),
Sergeyev&Markin(1995), Strongin&Sergeyev(2000)}).

Let us introduce now possible ways to fix  STEP 4 of the $GS$. At
this step, we select an interval where a new trial will be
executed. We   consider both the traditional rule used, for
example, in \cite{Piya(1972)} and \cite{Strongin&Sergeyev(2000)}
and a new one that we shall call  the \textit{local improvement}
technique. The   traditional way to choose an interval for the
next trial is the following. \vskip 12pt \noindent
{\textbf{STEP 4.1.} \\
\hspace*{1cm} Select the interval $(x_{t-1}, x_t)$  such that
\begin{equation} \label{mint}
R_t= \textrm{min} \{R_i \ : \ 2\leq i \leq k \}
\end{equation}
and $t$ is the minimal number satisfying (\ref{mint}). \vskip 12pt

This rule used together with STEP 2.1 gives us Piyavskii's
algorithm. In this case,  the new trial point $x^{k+1} \in
(x_{t-1}, x_t)$ is
 chosen in such a way that
 \[
    R_t =  \min \{ R_i : 2\leq i \leq k \} = c_t(x^{k+1}) = \min \{ C^k(x) : x \in [a,b] \}.
 \]
The new way to fix STEP 4 is introduced below. \vskip 12pt

\noindent {\textbf{STEP 4.2 (the  local improvement
technique).}} \\
\hspace*{1cm} $flag$ is a parameter initially equal to zero. \\
\hspace*{1cm} $imin$ is the index corresponding to the current
estimate of the minimal value \\
\hspace*{1cm} of the function, that is: $z_{imin}=f(x_{imin})\le f(x_i)$, $i=1,...,k$. \\
%\hspace*{1cm} $imin=argmin \{z_i \ : \ 1\leq i \leq k \}$\\
\hspace*{1cm} $z^k$ is the result of the last trial corresponding to a point $x_j$ in the line (\ref{order}),\\
\hspace*{1cm}  i.e., $x^k=x_j$. \\
\hspace*{1cm} IF (flag=1) THEN \\
\hspace*{2cm} IF $z^k<z_{imin}$ THEN $imin=j$.  \\
\hspace*{2cm} {\em Local improvement:} Alternate the choice of the
interval $(x_{t-1}, x_t)$  among\\
\hspace*{2cm} $t=imin+1$ and $t=imin$, if $imin=2,...,k-1,$ (if
$imin=1$ or $imin=k$\\
\hspace*{2cm}  take $t=2$ or $t=k$, respectively) in such a way
that for   $\delta > 0$  it follows
 \begin{equation} \label{delta}
 |x_t-x_{t-1}| > \delta. \hspace{.5cm}
 \end{equation}
\hspace*{1cm} ELSE (flag=0)\\
\hspace*{2cm} $ t= \textrm{argmin} \{R_i \ : \ 2\leq i \leq k \}$\\
\hspace*{1cm} ENDIF\\
\hspace*{1cm} flag=NOTFLAG(flag)

\vskip 12pt The motivation of the introduction of STEP 4.2
presented above is the following. In STEP 4.1, at each iteration,
we continue the search at an interval corresponding to the minimal
value of the characteristic $R_i$, $i=2,...,k$ (see (\ref{mint})).
This choice admits occurrence of such a situation where  the
search goes on for a certain finite (but possibly high) number of
iterations at subregions of the domain that are ``distant'' from
the best found approximation to the global solution  and only
successively concentrates trials at the interval containing a
global minimizer. However, very often it is of a crucial
importance to be able to find a good approximation of the global
minimum in the lowest number of iterations. Due to this reason, in
STEP 4.2 we take into account the rule (\ref{mint}) used in STEP
4.1 and related to the minimal characteristic, but we alternate it
with a new selection method that forces the algorithm to continue
the search in the part of the domain close to the best value of
the objective function found up to now. The parameter ``flag''
assuming values 0 or 1 allows us to alternate the two methods of
the selection.

More precisely, in STEP 4.2   we start by identifying  the index
$imin$ corresponding to the current minimum among the   found
values of the objective  function $f(x)$, and then we select the
interval $(x_{imin}, x_{imin+1})$ located on the right of the best
current point,   $x_{imin}$,  or the interval on the left of
$x_{imin}$, i.e., $(x_{imin-1}, x_{imin})$. STEP 4.2 keeps working
alternatively on the right and on the left of the current best
point $x_{imin}$ until a new trial point with value less than
$z_{imin}$ is found. The search   moves from the right  to the
left of the best found approximation trying to improve it.
However, since we are not sure that the found best approximation
$x_{imin}$ is really located in the neighborhood of a global
minimizer $x^*$, the local improvement is alternated in STEP 4.2
with the usual rule (\ref{mint}) providing so the global search of
new subregions possibly containing the global solution $x^*$. The
parameter $\delta$ defines the width of the intervals that can be
subdivided during the phase of the local improvement. Note that
the trial points produced during the phases of the local
improvement (obviously, there can be more than one phase in the
course of the search) are used during the further iterations of
the global search in the same way as the points produced during
the global phases.

Let us consider now possible combinations of   the different
choices of STEP 2 and STEP 4 allowing us to construct the
following six algorithms.

\vskip 12pt \noindent \bd
 \item
 {$ \bf PKC$}: $GS$ with STEP 2.1
and STEP 4.1 (\textbf{P}iyavskii's method with the a priori
\textbf{K}nown \textbf{C}onstant $L$).
 \item
  {$ \bf GE$}:
$GS$ with STEP 2.2 and STEP 4.1 (the method using the
\textbf{G}lobal \textbf{E}stimate of the Lipschitz constant $L$).
 \item
  {$\bf LT$}: $GS$
with STEP 2.3 and STEP 4.1 (the method executing the
\textbf{L}ocal \textbf{T}uning on the local Lipschitz constants).
 \item
  {$\bf PKC\_LI$}: $GS$
with STEP 2.1 and STEP 4.2 (\textbf{P}iyavskii's method with the a
priori \textbf{K}nown \textbf{C}onstant $L$ enriched by the
\textbf{L}ocal \textbf{I}mprovement technique).
 \item
  {$\bf
GE\_LI$}: $GS$ with STEP 2.2 and STEP 4.2 (the method using the
\textbf{G}lobal \textbf{E}stimate of $L$ enriched by the
\textbf{L}ocal \textbf{I}mprovement technique).
 \item
  {$\bf LT\_LI$}: $GS$ with STEP 2.3 and STEP 4.2
(the method executing the \textbf{L}ocal \textbf{T}uning on the
local Lipschitz constants enriched by the \textbf{L}ocal
\textbf{I}mprovement technique). \ed

\vskip 12pt Let us consider convergence properties  of the
introduced algorithms by studying an infinite trial sequence
$\{x^k\}$ generated by an algorithm belonging to the general
scheme $GS$ for solving problem  (\ref{p}), (\ref{fun}). We remind
that the algorithm $PKC$ is   Piyavskii's method and its
convergence properties have been   studied  in \cite{Piya(1972)}.
In order to start we need  the following definition.
\begin{defin}\label{def1}
Convergence to a point $x' \in (a,b)$ is said to be bilateral if
there exist two infinite subsequences of $\{x^k\}$ converging to
$x'$ one from the left, the other from the right.
\end{defin}

\begin{theorem}
\label{th1} Assume that the objective function $f(x)$ satisfies
the condition (\ref{fun}), and let $x'$ be any limit point of
$\{x^k\}$ generated by the $GE$ or by the $LT$ algorithm. Then the following
assertions hold:
\begin{itemize}
\item[1.] convergence to $x'$ is bilateral, if $x' \in (a,b)$;
\item[2.] $f(x^k) \geq f(x')$, for all trial points $x^k$, $k\geq 1$;
\item[3.] if there  exists another limit point $x''\neq x'$, then $f(x'')=f(x')$;
\item[4.] if the function $f(x)$ has a finite number of local minima in
$[a,b]$, then the point $x'$ is locally optimal;
\item[5.] (Sufficient conditions for  convergence to a global minimizer).
Let $x^*$ be a global minimizer of $f(x)$.
If there exists an iteration number $k^*$ such that for all $k>k^*$ the
 inequality
\begin{equation} \label{conv}
l_{j(k)} \geq  L_{j(k)}
\end{equation}
holds, where $L_{j(k)}$ is the Lipschitz constant for the interval
$[x_{j(k)-1}, x_{j(k)}]$ containing $x^*$, and $l_{j(k)}$ is its
estimate (see (\ref{hglob}) and (\ref{hloc})). Then the set of
limit points of the sequence $\{x^k\}$ coincides with the set of
global minimizers of the function $f(x)$.
\end{itemize}
\end{theorem}
{\em Proof.} The proofs of assertions 1--5 are analogous to the
proofs of Theorems 4.1--4.2 and Corollaries 4.1--4.4 from
\cite{Strongin&Sergeyev(2000)}. \hfill $\Box$

\begin{theorem}
\label{th2} Assertions 1--5 of Theorem \ref{th1} hold for the
algorithms $PKC\_LI$, $GE\_LI$, and $LT\_LI$  for a fixed finite
tolerance $\delta>0$ and $\varepsilon =0$, where $\delta$ is from
(\ref{delta}) and $\varepsilon$ is from (\ref{epsilon}).
\end{theorem}
{\em Proof.} Since $\delta>0$ and $\varepsilon =0$, the algorithms
$PKC\_LI$, $GE\_LI$, and $LT\_LI$  use the local improvement only
at the initial stage of the search   until the selected interval
$(x_{t-1},x_t)$ is greater than $\delta$. When $|x_t-x_{t-1}| \leq
\delta$ the interval cannot be divided by the local improvement
technique and   the selection criterion (\ref{mint}) is used.
Thus, since the one-dimensional search region   has a finite
length and $\delta$ is a fixed finite number, there exists a
finite iteration number $j$ such that at all iterations $k>j$ only
selection criterion (\ref{mint}) will be used.  As a result, at
the remaining part of the search, the methods $PKC\_LI$, $GE\_LI$,
and $LT\_LI$  behave themselves as the algorithms $PKC$, $GE$, and
$LT$, respectively. This consideration concludes the proof. \hfill
$\Box$

The next theorem ensures existence of the values of the parameter
$r$ satisfying condition (\ref{conv}) providing so that all global
minimizers of $f(x)$ will be located by the four proposed methods
that do not use the a priori known Lipschitz constant.

\begin{theorem}
\label{th3} For any function $f(x)$ satisfying (\ref{fun}) with $L
< \infty$ there exists a value $r^*$ such that for all $r>r^*$
condition (\ref{conv}) holds for the four algorithms $GE$, $LT$,
$GE\_LI$, and $LT\_LI$.
\end{theorem}
{\em Proof.} It follows from (\ref{hglob}), (\ref{hloc}), and the
finiteness of $\xi
> 0$ that  approximations of the Lipschitz constant $l_i$ in the
four methods  are always greater than zero. Since $L<\infty$ in
(\ref{fun}) and any positive value of the parameter $r$ can be
chosen in the scheme $GS$, it follows that there exists an $r^*$
such that condition (\ref{conv}) will be satisfied for all global
minimizers for $r>r^*$. This fact, due to Theorems~\ref{th1}
and~\ref{th2}, proves the theorem.\hfill$\Box$

\section{Six methods constructing smooth piece-wise quadratic auxiliary
functions for solving  problems with the Lipschitz first
derivative} \setcounter{equation}{0}

\begin{figure}[t]
\centering \
\centerline{\epsfxsize=12cm\mbox{\epsfbox{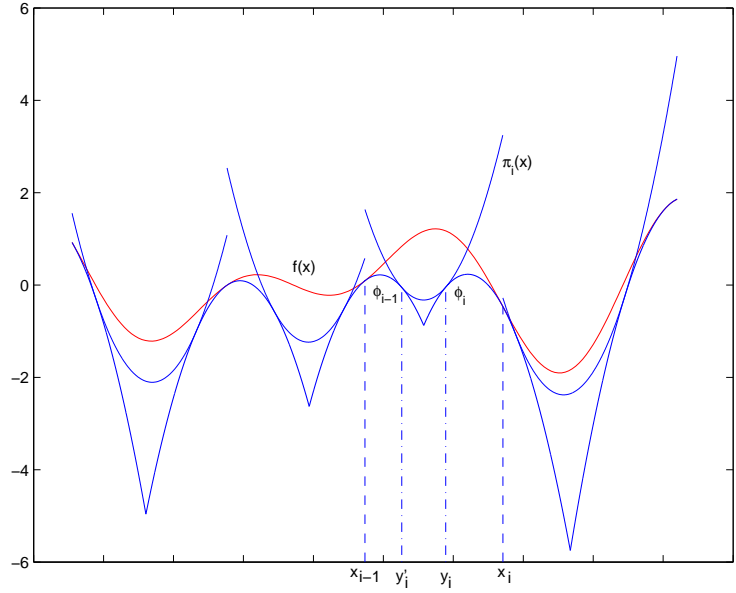}}}
\caption{\em Constructing smooth support functions by using
$\phi_{i-1}(x)$, $\pi_i(x)$, and $\phi_{i}(x)$} \label{fig.3-4-1}
\end{figure}

In this Section, we study the algorithms for solving problem
(\ref{p}) with the Lipschitz condition (\ref{der}) that holds for
the first derivative $f'(x)$ of the objective function $f(x)$. In
this Section, by the term {\em trial} we denote the evaluation of
both the function $f(x)$ and its first derivative $f'(x)$ at a
point $x$ that is called the {\em trial point}.

We consider the smooth support functions described in
\cite{Sergeyev(1998)}. This approach is based on the fact observed
in \cite{Sergeyev(1998)}, namely, that  at each interval
$[x_{i-1}, x_i]$  (see Fig.~\ref{fig.3-4-1})   the curvature of
the objective function $f(x)$ is determined by the Lipschitz
constant $M$ from (\ref{der}). In particular, over the interval
$(y'_i,y_i)$ it should be $f(x) \ge \pi_i(x)$ where
\begin{equation} \label{parab}
\pi_i(x)=0.5M x^2 + b_ix+c_i.
\end{equation}
This means that over the interval $(y'_i,y_i)$ both the objective
function $f(x)$ and the parabola $\pi_i(x)$ are strictly above the
Breiman-Cutler-Gergel's function $\Phi_i(x)$ from
(\ref{Cutlersupp})  where the unknowns $b_i, \ c_i, \ y'_i,$ and $
\ y_i$ can be   determined  following the considerations made in
\cite{Sergeyev(1998)}.

\begin{figure}
\centerline{\epsfxsize=12cm\mbox{\epsfbox{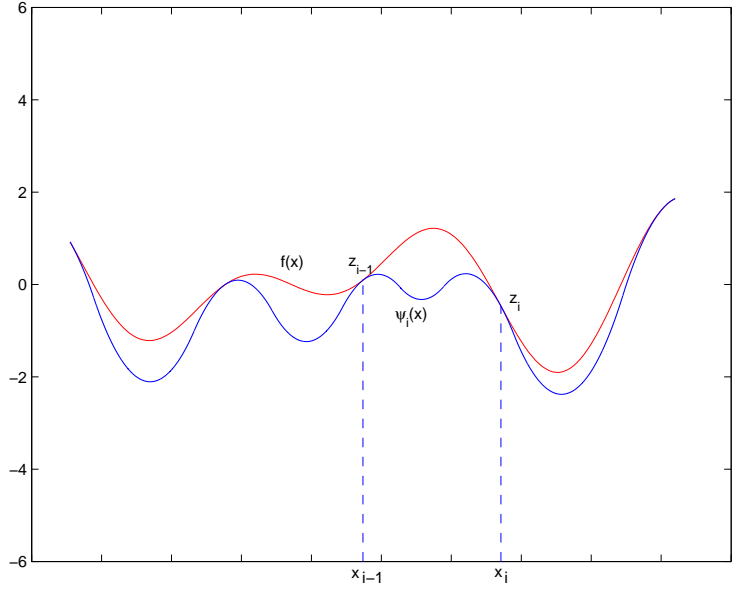}}}
 \caption{\em The resulting smooth support functions $\psi_i(x)$} \label{fig.3-4-2}
\end{figure}

These results from \cite{Sergeyev(1998)} allow us to construct the
following smooth support function $\psi_i(x)$ for $f(x)$ over
$[x_{i-1}, x_i]$:
\begin{equation} \label{suppsmooth}
\psi_i(x)=\left\{\begin{array}{lll}
               \phi_{i-1}(x), & x\in [x_{i-1},y'_i], \\
               \pi_i(x), & x\in [y'_i,y_i], \\
               \phi_i(x), & x\in [y_i, x_i]
        \end{array} \right.
\end{equation}
where there exists the first derivative $\psi'_i(x)$, $x\in
[x_{i-1}, x_i]$, and
$$ \psi_i(x) \leq f(x), \hspace{.5cm} x\in [x_{i-1}, x_i].  $$
This function is shown in Fig.~\ref{fig.3-4-2}. The points $y_i$,
$y'_i$ and the vertex $\bar x_i$ of the parabola $\pi_i(x)$ can be
found  (see \cite{Sergeyev(1998)} for the details) as follows:
\begin{equation}\label{yi}
y_i=\frac{x_i-x_{i-1}}4+\frac{z'_i-z'_{i-1}}{4M}+\frac{z_{i-1}-z_i+z'_ix_i
-z'_{i-1}x_{i-1}+0.5M(x_i^2-x_{i-1}^2)}{M(x_i-x_{i-1})+z'_i-z'_{i-1}},
\end{equation}
\begin{equation}\label{yprimoi}
y'_i=-\frac{x_i-x_{i-1}}4-\frac{z'_i-z'_{i-1}}{4M}+\frac{z_{i-1}-z_i
+z'_ix_i-z'_{i-1}x_{i-1}+0.5M(x_i^2-x_{i-1}^2)}{M(x_i-x_{i-1})+z'_i-z'_{i-1}},
\end{equation}
\begin{equation}\label{xbar}
\bar x_i=2y_i-\frac 1M z'_i-x_i,
\end{equation}
where $z_i=f(x_i)$ and $z'_i=f'(x_i)$.

In order to construct global optimization algorithms by applying
the same methodology used in the previous Section, for each
interval $[x_{i-1}, x_i]$ we should calculate its characteristic
$R_i$. For the smooth auxiliary functions $\psi_i(x)$ it can be
calculated as $R_i=\psi_i(p_i)$, where
$$ p_i= \arg \min \{\psi_i(x): \ x\in [x_{i-1}, x_i]\}. $$
Three different cases can take place.
\begin{itemize}
\item[i)]
 The first one is shown in Fig.~\ref{fig.3-4-2}. It
corresponds to the situation where conditions $\psi'_i(y'_i)<0$
and $\psi'_i(y_i)>0$ hold. In this case $$p_i= \arg \min
\{f(x_{i-1},\psi_i(\bar x_i),f(x_i)\}$$ and
\begin{equation} \label{car1}
R_i=\min \{f(x_{i-1}),\psi_i(\bar x_i),f(x_i)\}.
\end{equation}
\item[ii)]
 The second case is whenever $\psi'_i(y'_i)\geq 0$ and $\psi'_i(y_i)>0$. In this
 situation,
  we have (see \cite{Sergeyev(1998)}) that
\begin{equation} \label{car2}
R_i=\min \{f(x_{i-1}),f(x_i)\}.
\end{equation}
\item[iii)]
 The third case is when $\psi'_i(y'_i)< 0$ and $\psi'_i(y_i)\leq 0$. It can be
 considered by a complete analogy with the previous one.
\end{itemize}

We are ready now to introduce the general scheme for the methods
working with smooth piece-wise quadratic auxiliary functions. As
it has been done in the previous Section, six different algorithms
will be then constructed by specifying STEP 2 and STEP 4 of the
general scheme.

 \vskip 12pt {\bf General
Scheme   describing  algorithms working with the first Derivatives
and constructing smooth piece-wise quadratic auxiliary functions
($GS\_D$).}
 \bd \item[STEP 0.] The
first two trials are performed at the points $x^1=a$ and $x^2=b$.
The point $x^{k+1}$, $k\geq 2$, of the current (k+1)-th iteration
is chosen as follows. \item[STEP 1.] Renumber the trial points
$x^1, \ x^2, \ \dots, \ x^k$ of the previous iterations by
subscripts so that
 \begin{equation} \label{orderiv}
 a=x_1 < \dots < x_k =b.
\end{equation}
\item[STEP 2.] Compute in a certain way the values $m_i$  being
estimates of the Lipschitz constants of $f'(x)$ over the intervals
$[x_{i-1}, x_i]$, $i=2,...k.$ The way to calculate the values
$m_i$ will be specified in each concrete algorithm described
below.
 \item[STEP 3.]
Initiate the index sets $I=\emptyset$, $Y'=\emptyset$, and
$Y=\emptyset$. Set the index of the current interval $i=2$ and go
to STEP 3.1.
 \bd \item[STEP 3.1.] If for the current interval
$[x_{i-1}, x_i]$ the following inequality
\begin{equation} \label{sigpar}
 \pi'_i(y'_i)\cdot\pi'_i(y_i) \geq 0
\end{equation}
does not hold (where $\pi'(x)$ is the derivative of the parabola
(\ref{parab})) then go to STEP 3.2. Otherwise go to STEP 3.3.
\item[STEP 3.2.] Calculate for the interval $[x_{i-1}, x_i]$ its
characteristic $R_i$ using  (\ref{car1}). Include $i$ in $I$ and
go to STEP 3.4. \item[STEP 3.3.] Calculate for the interval
$[x_{i-1}, x_i]$ its characteristic $R_i$ using  (\ref{car2}). If
$$ f(x_{i-1})<f(x_i) $$
then include the index $i$ in the set $Y'$ and go to STEP 3.4.
Otherwise include $i$ in the set $Y$ and go to STEP 3.4.
\item[STEP 3.4.] If $i<k$, set $i=i+1$ and go to STEP 3.1.
Otherwise go to STEP~4. \ed
 \item[STEP 4.] Find the interval
$(x_{t-1}, x_t)$ for the next possible trial. The way to do it
will be specified in each concrete algorithm described below.
\item[STEP 5.] If
\begin{equation} \label{epsi}
 |x_t - x_{t-1} | > \varepsilon,
\end{equation}
where $\varepsilon > 0$ is a given search accuracy, then execute the next trial at the point
\begin{equation} \label{pointder}
x^{k+1} = \left\{\begin{array}{lll}
               y'_t,  & \mbox{ if } t\in Y', \\
               \bar x_t,  & \mbox{ if } t\in I, \\
               y_t,  & \mbox{ if } t\in Y,
        \end{array} \right.
\end{equation}
and go to STEP 1. Otherwise, take as an estimate of the global
minimum $f^*$ from (\ref{p}) the value
$$ f^*_k = \min \{z_i \ : \ 1\leq i \leq k \}, $$
and a   point
$$ x^*_k = \arg \min \{z_i \ : \ 1\leq i \leq k \}, $$
as an estimate of the global minimizer $x^*$, after executing
these operations STOP.
  \ed

\vskip 12pt Let us make just two comments upon the introduced
scheme $GS\_D$. First,   in STEPS 3.1--3.4 the characteristics
$R_i$, $i=2,...,k$, are calculated by taking into account the
different cases i -- iii of the location of the point $p_i$
described above. Second, note that the index sets $I$, $Y$, and
$Y'$ have been introduced in order to calculate the new trial
point $x^{k+1}$ in STEP 5. In fact, the vertex $\bar x_i$ of the
$i-th$ parabola, $i=2,...,k$, can  be  outside the interior of the
interval $[x_{i-1}, x_i]$. It can happen that $\bar x_i \notin
[x_{i-1}, x_i]$ whenever $\psi'_i(y'_i)\geq 0$ and
$\psi'_i(y_i)>0$ (or $\psi'_i(y'_i)< 0$ and $\psi'_i(y_i)\leq 0$),
and so the point $y'_i$ (or $y_i$) is selected as new trial point
$x^{k+1}$.

 \vskip 12pt Let us show now how it is possible to specify
STEP 2 and STEP 4 of the scheme $GS\_ D$. As it has been done in
the previous Section for the scheme $GS$, we first describe three
different choices of the values $m_i$ that should be done at
STEP~2 and then consider  two selection rules  that can be used to
fix STEP~4 for choosing the point $x^{k+1}$. The first possible
way to assign values to $m_i$ is the following: \vskip12pt
\noindent
{\textbf{STEP 2.1} } \\
\hspace*{1cm} Set
\begin{equation} \label{hnotader}
 m_i= M, \hspace{1cm} i=2,...,k,
\end{equation}
\hspace*{1cm} where $M$ is from (\ref{der}). \vskip12pt

In this case, the exact value of the a priori given Lipschitz
constant for the first derivative $f'(x)$ is used. As a result,
the auxiliary functions $\psi_i(x)$ from (\ref{suppsmooth}) are
support functions for $f(x)$ over the intervals $[x_{i-1}, x_i]$,
$i=2,...,k$. Since it is difficult to know the exact value $M$ in
practice, the choices made  in the following STEPS 2.2 and 2.3 (as
it was for the methods working with Lipschitz objective functions)
describe how to estimate dynamically the global constant $M$ (STEP
2.2) and the local constants related to each interval $[x_{i-1},
x_i]$, $i=2,...,k$ (STEP  2.3). \vskip12pt \noindent
{\textbf{STEP 2.2}  } \\
\hspace*{1cm} Set
\begin{equation} \label{hglobder}
 m_i= r \max \{\xi, H^k\}, \hspace{1cm} i=2,...,k,
\end{equation}
\hspace*{1cm} where  $\xi >0$ reflects the supposition that $f'(x)$ is not
constant over the interval \\
\hspace*{1cm} $[a,b]$ and $r>1$ has the same sense as in the
STEP~2.2 of the scheme $GS$. The\\
\hspace*{1cm} value $H^k$  is computed as
\begin{equation} \label{Hmaxder}
H^k= \max \{v_i: i=2,...,k\},
\end{equation}
\hspace*{1cm} where
\begin{equation} \label{vi}
 v_i=\frac{|2(z_{i-1}-z_i)+(z'_{i-1}+z'_i)(x_i - x_{i-1})|+d_i}{(x_i - x_{i-1})^2}
\end{equation}
\hspace*{1cm} and
\begin{equation} \label{di}
 d_i=\sqrt{|2(z_{i-1}-z_i)+(z'_{i-1}+z'_i)(x_i - x_{i-1})|^2+(z'_i-z'_{i-1})^2(x_i - x_{i-1})^2}.
\end{equation}

If an algorithm uses the exact value $M$ of the Lipschitz constant
(see STEP 2.1 above) then it is ensured by construction that the
points $y'_i$, $y_i$ from (\ref{yprimoi}) and (\ref{yi}) belong to
the interval $[x_{i-1}, x_i]$. In the case, when an estimate $m_i$
of $M$   is used, it can happen that, if the value $m_i$ is
underestimated,  the points $y'_i$ and $y_i$ can be obtained
outside the interval $[x_{i-1}, x_i]$ that would lead to an error
in the work of the algorithm using such an underestimate. It has
been proved in \cite{Sergeyev(1998)} that the choice
(\ref{hglobder})--(\ref{di}) makes this unpleasant situation
impossible. More precisely, the following theorem holds.
\begin{theorem}\label{teo1}
If the values $m_i$ in $GS\_D$ are determined by formulae
(\ref{hglobder})-(\ref{di})  then the points $y'_i$, $y_i$ from
  (\ref{yi}), (\ref{yprimoi})  belong to the interval
$[x_{i-1}, x_i]$ and the following estimates take place:
$$ y'_i-x_{i-1} \geq \frac{(r-1)^2}{4r(r+1)}(x_i - x_{i-1}), $$
$$ x_i-y_i \geq \frac{(r-1)^2}{4r(r+1)}(x_i - x_{i-1}). $$
\end{theorem}

Let us introduce now STEP  2.3 that shows how the local tuning
technique works in the situation where the first derivative of the
objective function can be calculated.
 \vskip 12pt \noindent
{\textbf{STEP 2.3}  } \\
\hspace*{1cm} Set
\begin{equation} \label{hlocder}
 m_i = r \max \{\lambda_i, \gamma_i, \xi \},
\end{equation}
\hspace*{1cm} where $r>1$ and $\xi>0$ have the same sense as before, and
\begin{equation} \label{landader}
 \lambda_i= \max\{v_{i-1}, v_i, v_{i+1}\}, \ \ i=3,...,k-1,
\end{equation}
\hspace*{1cm} where the values $v_i$ are calculated following
(\ref{vi}), and when $i=2$ and $i=k$ \\ \hspace*{1cm} we consider
only $v_2, v_3$, and $v_{k-1}, v_k$, respectively.  The value
$\gamma_i$ is computed \\ \hspace*{1cm} as follows
\begin{equation} \label{gammader}
\gamma_i = H^k \frac{(x_i-x_{i-1})}{ X^{max}},
\end{equation}
\hspace*{1cm} where $H^k$ is from (\ref{Hmaxder}) and
\[
 X^{max} = \max \{(x_i-x_{i-1}), \ 1=2,...,k \}.
\]

\vskip 12pt   As it was in STEP 2.3  of the scheme $GS$ from the
previous Section, the local tuning technique balances the local
and the global information to get the estimates $m_i$ on the basis
of the local and the global estimates $\lambda_i$ and $\gamma_i$.
Note also that the fact that $y'_i$ and $y_i$ belong to the
interval $[x_{i-1}, x_i]$ can be proved by a complete analogy with
Theorem \ref{teo1} above.

 \vskip 12pt Let us consider now  STEP 4 of the scheme $GS\_D$. At this step, we should select
an interval $[x_{t-1}, x_t]$ containing the next trial point
$x^{k+1}$. As we have already done in   Section 2, we consider two
strategies: the rule   selecting the interval corresponding to the
minimal characteristic $R_t$ and the local improvement technique.
Thus, STEP 4.1 and STEP 4.2 of the scheme $GS\_D$
  correspond exactly to   STEP 4.1   and STEP 4.2 of the scheme
$GS$  from Section 2. The obvious difference consists of the fact
that characteristics $R_i$, $i=2,...,k$ are calculated with
respect to STEPS 3.1--3.4 of the scheme $GS\_D$.

 \vskip 12pt Thus,   by specifying STEP 2 and STEP 4 we obtain from the general scheme
$GS\_D$ the following six algorithms:
 \bd
\item {$ \bf DKC$}: $GS\_D$ with STEP 2.1 and STEP 4.1 (the method
using the first \textbf{D}erivatives and the a priori
\textbf{K}nown Lipschitz \textbf{C}onstant $M$).
 \item
{$ \bf DGE$}: $GS\_D$ with STEP 2.2 and STEP 4.1 (the method using
the first \textbf{D}erivatives and the \textbf{G}lobal
\textbf{E}stimate of the Lipschitz constant $M$).
 \item
{$\bf DLT$}: $GS\_D$ with STEP 2.3 and STEP 4.1 (the method using
the first \textbf{D}erivatives and the \textbf{L}ocal
\textbf{T}uning).
 \item
{$\bf DKC\_LI$}: $GS\_D$ with STEP 2.1 and STEP 4.2 (the method
using the first \textbf{D}erivatives,   the a priori
\textbf{K}nown Lipschitz \textbf{C}onstant $M$, and the
\textbf{L}ocal \textbf{I}mprovement technique).
 \item
{$\bf DGE\_LI$}: $GS\_D$ with STEP 2.2 and STEP 4.2 (the method
using the first \textbf{D}erivatives, the \textbf{G}lobal
\textbf{E}stimate of the Lipschitz constant $M$, and the
\textbf{L}ocal \textbf{I}mprovement technique).
 \item
{$\bf DLT\_LI$}: $GS\_D$ with STEP 2.3 and STEP 4.2 (the method
using the first \textbf{D}erivatives, the \textbf{L}ocal
\textbf{T}uning, and the \textbf{L}ocal \textbf{I}mprovement
technique).
 \ed

 Let us consider now   infinite trial sequences
$\{x^k\}$ generated by methods belonging to the general scheme
$GS\_D$ and study convergence properties of the six algorithms
introduced above.
\begin{theorem}
\label{th4} Assume that the objective function $f(x)$ satisfies
  condition (\ref{der}), and let $x'$ ($x'\neq a$, $x'\neq b$) be any limit point of
$\{x^k\}$ generated by either by the method $DKC$ or the $DGE$ or
the $DLT$. If the values $m_i$, $i=2,...,k$,  are bounded as below
\begin{equation} \label{convm}
v_i \leq m_i     < \infty,
\end{equation}
where $v_i$ is from (\ref{vi}), then the following assertions hold:
\begin{itemize}
\item[1.] convergence to $x'$ is bilateral, if $x' \in (a,b)$;
\item[2.] $f(x^k) \geq f(x')$, for all trial points $x^k$, $k\geq
1$; \item[3.] if there  exists another limit point $x''\neq x'$,
then $f(x'')=f(x')$;
 \item[4.] if the function $f(x)$ has a finite
number of local minima in $[a,b]$, then the point $x'$ is locally
optimal;
 \item[5.] (Sufficient conditions for  convergence to a
global minimizer). Let $x^*$ be a global minimizer of $f(x)$ and
$[x_{j(k)-1}, x_{j(k)}]$ be an interval containing this point
during the course of the $k$-th iteration of one of the algorithms
$DKC$, $DGE$, or $DLT$. If there exists an iteration number $k^*$
such that for all $k>k^*$ the inequality
\begin{equation} \label{convder}
 M_{j(k)} \leq m_{j(k)}  < \infty
\end{equation}
takes places for $[x_{j(k)-1}, x_{j(k)}]$ and  (\ref{convm}) for
all the other intervals, then the set of limit  points of the
sequence $\{x^k\}$ coincides with the set of global minimizers of
the function $f(x)$.
\end{itemize}
\end{theorem}
{\em Proof.} The proofs of assertions 1--5 are analogous to the
proofs of Theorems 5.1--5.5 and Corollaries 5.1--5.6 from
\cite{Sergeyev(1998)}. \hfill $\Box$

The fulfillment of the sufficient conditions for  convergence to a
global minimizer, i.e.,  (\ref{convder}), are evident for the
algorithm $DKC$. For the methods $DGE$ and $DLT$, its fulfillment
depends on the choice of the reliability parameter $r$. A theorem
similar to the theorem~\ref{th3} can be proved for them by a
complete analogy. However, there exist  particular cases where the
objective function $f(x)$ is such that its structure ensures that
(\ref{convder}) holds. In the following theorem, sufficient
conditions providing the fulfillment of (\ref{convder}) for the
methods $DGE$ and $DLT$ are established for a particular class of
objective functions. The theorem states that if $f(x)$ is quadratic
in a neighborhood $I(x^*)$ of the global minimizer $x^*$, then to
ensure the global convergence it is sufficient that   the methods
will place one trial point on the left from $x^*$ and one trial
point on the right from $x^*$.
\begin{theorem}
\label{th5} If the objective function $f(x)$ is such that there
exists  a neighborhood $I(x^*)$ of a global minimizer $x^*$ where
\begin{equation} \label{convr}
f(x)= 0.5 M x^2 + qx + n,
\end{equation}
where $q$ and $n$ are finite constants and $M$ is from (\ref{der})
and trials have been executed at points $x^-$, $x^+$ $\in I(x^*)$,
then condition (\ref{convder}) holds for algorithms $DGE$ and $DLT$
and $x^*$ is a limit point of the trial sequences generated by these
methods if (\ref{convm}) is fulfilled for all the other intervals.
\end{theorem}
{\em Proof.} The proof is  analogous to the proof of Theorem 5.6
from \cite{Sergeyev(1998)}. \hfill $\Box$ \vskip 12pt \noindent
\begin{theorem}
\label{th6} Assertions 1--5 of Theorem \ref{th4} hold for the
algorithms $DKC\_LI$, $DGE\_LI$, and $DLT\_LI$  for a fixed finite
tolerance $\delta>0$ and $\varepsilon =0$, where $\delta$ is from
(\ref{delta}) and $\varepsilon$ is from (\ref{epsilon}).
\end{theorem}
{\em Proof.} The proof is  analogous to the proof of Theorem
\ref{th2} from Section 2. \hfill $\Box$

\section{Numerical experiments}
\setcounter{equation}{0}

In this section, we present numerical
results executed on 120 functions taken from the literature to
compare the performance of the six algorithms described in Section
2 and the six algorithms from Section 3.

Two series of experiments have been done. In both of them the choice
of the reliability parameter $r$ has been done with the step 0.1
starting from $r=1.1$, i.e., $r=1.1, 1.2$, etc. in order to ensure
convergence to the global solution for all the functions taken into
consideration in each series. It is well known  (see detailed
discussions on the choice of $r$ and its influence on the speed of
Lipschitz global optimization methods in
\cite{Pinter(1996),Sergeyev&Kvasov(2008),Strongin&Sergeyev(2000)})
that in general, for higher values of $r$ methods of this kind are
more reliable but slower. It can be seen from the results of
experiments (see Tables~\ref{table1} -- \ref{table6}) that the
tested methods were able to find the global solution already for
very low values of $r$. Then, since there is no sense to make a
local improvement with the accuracy $\delta$ that is higher than the
final required accuracy $\varepsilon$, in all the algorithms using
the local improvement technique the accuracy $\delta$ from
(\ref{delta}) has been fixed $\delta=\varepsilon$. Finally, the
technical parameter $\xi$ (used only when at the initial iterations
a method executes trials at the points with equal values) has been
fixed to $\xi = 10^{-8}$ for all the methods using it.

In the first series of experiments, a set of 20 functions described
in \cite{Hansen&Jaumard(1992)} has been considered. In Tables
\ref{table1} and \ref{table2}, we present numerical results for the
six methods proposed to work with the problem (\ref{p}),
(\ref{fun}). In particular, Table \ref{table1} contains the numbers
of trials executed by the algorithms with the accuracy $\varepsilon
= 10^{-4}(b-a)$, where $\varepsilon$ is from (\ref{epsilon}). Table
\ref{table2} presents the results for $\varepsilon = 10^{-6}(b-a)$.
The parameter  $r=1.1$ was sufficient for the algorithms $GE$, $LT$,
and $GE\_LI$, $LT\_LI$ while the exact values of the Lipschitz
constant of the functions $f(x)$  have been used in the methods
$PKC$ and $PKC\_LI$.

\begin{table}[t]
\centering \footnotesize
\begin{tabular}{|c|cccccc|}
\hline \hline {\em Problem}  & {\em PKC} & {\em GE} & {\em LT} &
{\em  PKC\_LI } & {\em  GE\_LI } &  {\em  LT\_LI }
\\ \hline
1   & 149 & 158 & 37 & 37 & 35 & 35\\
2   & 155 & 127 & 36 & 33 & 35 & 35\\
3   & 195 & 203 & 145 & 67 & 25 & 41\\
4   & 413 & 322 & 45 & 39 & 39 & 37\\
5   & 151 & 142 & 46 & 151 & 145 & 53\\
6   & 129 & 90 & 84 & 39 & 41 & 41\\
7   & 153 & 140 & 41 & 41 & 33 & 35\\
8   & 185 & 184 & 126 & 55 & 41 & 29\\
9   & 119 & 132 & 44 & 37 & 37 & 35\\
10  & 203 & 180 & 43 & 43 & 37 & 39\\
11  & 373 & 428 & 74 & 47 & 43 & 37\\
12  & 327 & 99 & 71 & 45 & 33 & 35\\
13  & 993 & 536 & 73 & 993 & 536 & 75\\
14  & 145 & 108 & 43 & 39 & 25 & 27\\
15  & 629 & 550 & 62 & 41 & 37 & 37\\
16  & 497 & 588 & 79 & 41 & 43 & 41\\
17  & 549 & 422 & 100 & 43 & 79 & 81\\
18  & 303 & 257 & 44 & 41 & 39 & 37\\
19  & 131 & 117 & 39 & 39 & 31 & 33\\
20  & 493 & 70 & 70 & 41 & 37 & 33\\ \hline  {\em Average} &
314.60 & 242.40 & 65.10 & 95.60 & 68.55 & 40.80
\\ \hline\hline
\end{tabular}
\caption {\label{table1} \em  Results of numerical experiments
executed on 20 test problems from \cite{Hansen&Jaumard(1992)}  by
the six methods belonging to the scheme $GS$;  the accuracy
$\varepsilon=10^{-4}(b-a)$, $r=1.1$}
\end{table}

\begin{table}[t]
\centering \footnotesize
\begin{tabular}{|c|cccccc|}
\hline \hline {\em Problem}  & {\em PKC} & {\em GE} & {\em LT} &
{\em  PKC\_LI } & {\em  GE\_LI } &  {\em  LT\_LI }
\\ \hline
1   & 1681 & 1242 & 60 & 55 & 55 & 57\\
2   & 1285 & 1439 & 58 & 53 & 61 & 57\\
3   & 1515 & 1496 & 213 & 89 & 51 & 61\\
4   & 4711 & 3708 & 66 & 63 & 63 & 59\\
5   & 1065 & 1028 & 67 & 59 & 65 & 74\\
6   & 1129 & 761 & 81 & 63 & 65 & 61\\
7   & 1599 & 1362 & 64 & 65 & 55 & 59\\
8   & 1641 & 1444 & 194 & 81 & 67 & 49\\
9   & 1315 & 1386 & 64 & 61 & 59 & 57\\
10  & 1625 & 1384 & 65 & 59 & 63 & 57\\
11  & 4105 & 3438 & 122 & 71 & 63 & 61\\
12  & 3351 & 1167 & 114 & 67 & 57 & 55\\
13  & 8057 & 6146 & 116 & 8057 & 6146 & 119\\
14  & 1023 & 1045 & 66 & 57 & 49 & 49\\
15  & 7115 & 4961 & 103 & 65 & 61 & 59\\
16  & 4003 & 6894 & 129 & 63 & 65 & 63\\
17  & 5877 & 4466 & 143 & 69 & 103 & 103\\
18  & 3389 & 2085 & 67 & 65 & 61 & 57\\
19  & 1417 & 1329 & 60 & 61 & 57 & 53\\
20  & 2483 & 654 & 66 & 61 & 61 & 53\\ \hline  {\em Average} &
2919.30 & 2371.75 & 95.90 & 464.20 & 366.35 & 63.15
\\ \hline\hline
\end{tabular}
\caption {\label{table2} \em  Results of numerical experiments
executed on 20 test problems from \cite{Hansen&Jaumard(1992)} by
the six methods belonging to the scheme $GS$; the accuracy
$\varepsilon=10^{-6}(b-a)$,  $r=1.1$}
\end{table}

\begin{table}[t]
\centering \footnotesize
\begin{tabular}{|c|cccccc|}
\hline \hline {\em Problem}  & {\em DKC} & {\em DGE} & {\em DLT} &
{\em DKC\_LI} & {\em DGE\_LI} &  {\em DLT\_LI}
\\ \hline
1   & 15 & 16 & 14 & 15 & 16 & 15\\
2   & 10 & 12 & 12 & 10 & 13 & 13\\
3   & 48 & 58 & 56 & 57 & 27 & 25\\
4   & 14 & 14 & 11 & 14 & 14 & 11\\
5   & 17 & 16 & 15 & 17 & 20 & 18\\
6   & 22 & 24 & 22 & 21 & 25 & 21\\
7   & 10 & 12 & 11 & 11 & 13 & 11\\
8   & 44 & 52 & 50 & 45 & 25 & 25\\
9   & 10 & 13 & 12 & 10 & 13 & 13\\
10  & 9 & 12 & 12 & 9 & 12 & 12\\
11  & 24 & 26 & 22 & 26 & 25 & 21\\
12  & 19 & 21 & 21 & 19 & 23 & 21\\
13  & 197 & 63 & 25 & 37 & 45 & 26\\
14  & 13 & 17 & 17 & 14 & 20 & 20\\
15  & 55 & 43 & 17 & 43 & 33 & 18\\
16  & 66 & 54 & 26 & 49 & 45 & 26\\
17  & 51 & 45 & 32 & 33 & 39 & 27\\
18  & 5 & 9 & 9 & 5 & 9 & 9\\
19  & 11 & 11 & 11 & 11 & 11 & 11\\
20  & 22 & 24 & 25 & 19 & 23 & 25\\ \hline  {\em Average} & 33.10
& 27.10 & 21.00 & 23.25 & 22.55 & 18.40
\\ \hline\hline
\end{tabular}
\caption {\label{table3} \em Results of numerical experiments
executed on 20 test problems from \cite{Hansen&Jaumard(1992)}  by
the six methods belonging to the scheme $GS\_D$;  the accuracy
$\varepsilon=10^{-4}(b-a)$,   $r=1.2$}
\end{table}
\begin{table}[t]
\centering \footnotesize
\begin{tabular}{|c|cccccc|}
\hline \hline {\em Problem}  & {\em DKC} & {\em DGE} & {\em DLT} &
{\em DKC\_LI} & {\em DGE\_LI} &  {\em DLT\_LI}
\\ \hline
1   & 22 & 22 & 16 & 22 & 22 & 17\\
2   & 12 & 16 & 16 & 12 & 17 & 17\\
3   & 57 & 69 & 66 & 63 & 37 & 35\\
4   & 19 & 19 & 16 & 19 & 20 & 17\\
5   & 21 & 19 & 18 & 21 & 24 & 22\\
6   & 26 & 27 & 26 & 27 & 28 & 28\\
7   & 12 & 16 & 15 & 13 & 18 & 16\\
8   & 48 & 61 & 60 & 49 & 33 & 33\\
9   & 12 & 16 & 15 & 12 & 16 & 16\\
10  & 11 & 16 & 16 & 11 & 17 & 17\\
11  & 37 & 37 & 29 & 37 & 35 & 29\\
12  & 27 & 31 & 28 & 29 & 31 & 27\\
13  & 308 & 93 & 29 & 55 & 61 & 30\\
14  & 16 & 21 & 21 & 17 & 25 & 25\\
15  & 87 & 66 & 21 & 67 & 51 & 22\\
16  & 101 & 81 & 28 & 65 & 63 & 28\\
17  & 73 & 66 & 38 & 51 & 61 & 39\\
18  & 5 & 14 & 14 & 5 & 14 & 14\\
19  & 13 & 15 & 15 & 13 & 15 & 15\\
20  & 24 & 27 & 27 & 25 & 28 & 28\\ \hline
{\em Average} & 46.55 &
36.60 & 25.70 & 30.65 & 30.80 & 23.75
\\ \hline\hline
\end{tabular}
\caption {\label{table4} \em Results of numerical experiments
executed on 20 test problems from \cite{Hansen&Jaumard(1992)}  by
the six methods belonging to the scheme $GS\_D$;  the accuracy
$\varepsilon=10^{-6}(b-a)$,   $r=1.2$}
\end{table}

In Tables \ref{table3} and \ref{table4}, we present numerical
results for the six methods proposed to work with the problem
(\ref{p}), (\ref{der}). We have considered the same two accuracies
in Tables \ref{table1} and \ref{table2}: $\varepsilon=10^{-4}(b-a)$
for the experiments shown in Table \ref{table3} and
$\varepsilon=10^{-6}(b-a)$ for the results presented in Table
\ref{table4}. The reliability parameter $r$ has been taken equal to
$1.2$.

All the global minima have been found by all the methods in all
the experiments presented in Tables \ref{table1}--\ref{table4}. In
the last rows of  Tables \ref{table1} -- \ref{table4}, the average
values of the numbers of trials points generated by the algorithms
are given. The first (quite obvious) observation that can be made
with respect to the performed experiments consists of the fact
that the methods using derivatives are faster than the methods
that do not use this information (compare results in Tables
\ref{table1} and \ref{table2} with the results in Tables
\ref{table3} and \ref{table4}, respectively).

Then, it can be seen from Tables \ref{table1} and \ref{table2}
that both accelerating techniques, the local tuning and the local
improvement, allow us to speed up the search significantly when we
work with the methods belonging to the scheme $GS$. With respect
to the local tuning we can see that the method $LT$ is faster than
the algorithms $PKC$ and $GE$. Analogously, the $LT\_LI$ is faster
than the methods $PKC\_LI$ and $GE\_LI$. The introduction of the
local improvement also was very successful. In fact, the
algorithms $PKC\_LI$, $GE\_LI$, and $LT\_LI$ work significantly
faster than the  methods $PKC$, $GE$, and $LT$, respectively.

The effect of the local tuning and local improvement techniques is
very well marked in the case of the problem (\ref{p}), (\ref{fun}),
i.e., when the first derivative of the objective function is not
available. In the case of methods proposed to solve the problem
(\ref{p}), (\ref{der}), where the derivative of the objective
function can be used to construct algorithms, the effect  of the
introduction of the two acceleration techniques is always present
but is not so strong (see Tables \ref{table3} and \ref{table4}).
This happens because the smooth auxiliary functions constructed by
all the methods belonging to the scheme $GS\_D$ are much better than
the piece-wise linear functions build by the methods belonging to
the class $GS$. The smooth auxiliary functions are very close to the
objective function $f(x)$ providing so already in the case of the
slowest $DKC$ algorithm a very good speed and, as a result, leaving
less space for a possible acceleration that can be obtained thanks
to applying the local tuning and/or the local improvement
techniques. However, also in this case, the algorithm $DLT\_LI$ is
two times faster than the method $DKC$ (see Tables \ref{table3} and
\ref{table4}). Finally, it can be clearly seen from Tables
\ref{table1} -- \ref{table4} that the acceleration effects produced
by both techniques are more pronounced when the accuracy of the
search increases. This effect takes place for the methods belonging
to both schemes, $GS$ and  $GS\_D$.

\begin{figure}[t]
\centering
\centerline{\epsfxsize=10.5cm\mbox{\epsfbox{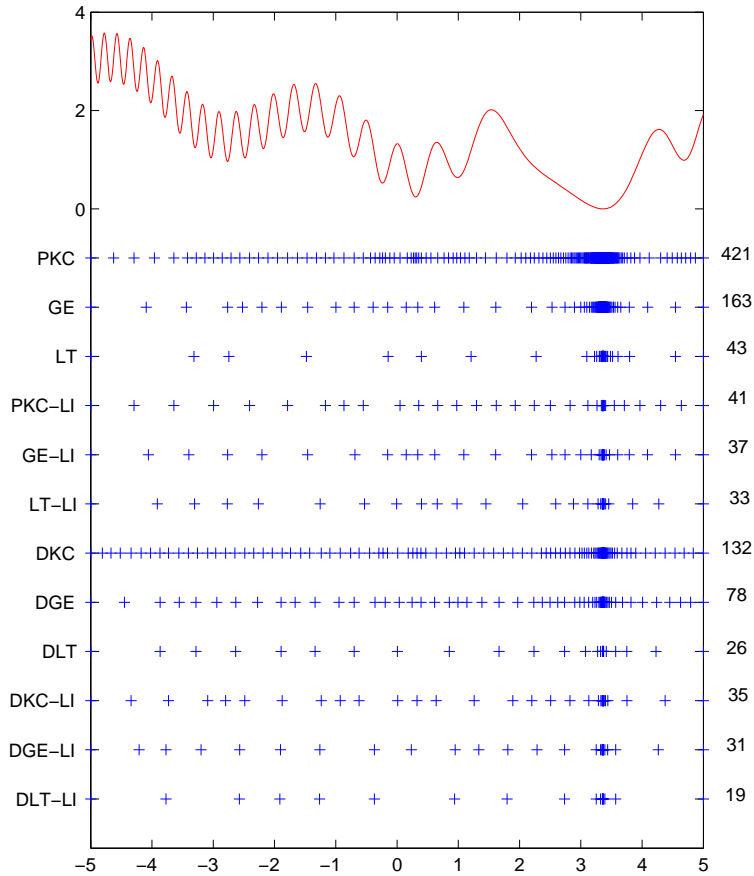}}}
\caption{\em Graph of  the function number 38 from (\ref{fpinter})
and trial points generated by the 12 methods while minimizing this
function.} \label{fig.5}
\end{figure}
\begin{table}[h]
\centering \footnotesize
\begin{tabular}{lllll}
\hline {\em Method} & {\em $r$} & {\em $\varepsilon=10^{-4}$} &
{\em $r$} & {\em $\varepsilon=10^{-6}$}
\\ \hline \\
PKC   &  & 400.54 &  & 2928.48 \\
GE & 1.1 & 167.63 & 1.1 & 1562.27 \\
LT & 1.1 & 47.28 & 1.1 & 70.21 \\
PKC\_LI &  & 44.82 &  & 65.70 \\
GE\_LI & 1.1 & 40.22 & 1.2 & 62.96 \\
LT\_LI & 1.3* & 38.88 & 1.2 & 60.04
\\ \hline
\end{tabular}
\caption {\label{table5} \em Average number of trial points
generated by   the six methods belonging to the scheme $GS$ on 100
test functions from \cite{Pinter(2002)} with the accuracies
$\varepsilon=10^{-4}$ and $\varepsilon=10^{-6}$}
%\\ (*) The value r=1.3 is used
%for 99 functions, for the function no.32 we have considered the value r=1.4}
\end{table}

\begin{table}[h]
\centering \footnotesize
\begin{tabular}{lllll}
\hline {\em Method} & {\em $r$} & {\em $\varepsilon=10^{-4}$} &
{\em $r$} & {\em $\varepsilon=10^{-6}$}
\\ \hline \\
DKC   &  & 125.85 &  & 170.65 \\
DGE & 1.1 & 87.53 & 1.1 & 121.01 \\
DLT & 1.1 & 49.00 & 1.1 & 53.53 \\
DKC\_LI &  & 43.72 &  & 62.88 \\
DGE\_LI & 1.1 & 38.46 & 1.1 & 58.61 \\
DLT\_LI & 1.1 & 28.50 & 1.1 & 40.57
\\ \hline
\end{tabular}
\caption {\label{table6} \em Average number of trial points
generated by   the six methods belonging to the scheme $GS\_D$ on
100 test functions from \cite{Pinter(2002)} with the accuracies
$\varepsilon=10^{-4}$ and $\varepsilon=10^{-6}$}
\end{table}

 In the second series of
experiments, a class of 100 one-dimensional randomized test
functions from \cite{Pinter(2002)} has been taken. Each function
$f_j(x)$, $j=1,...,100$, of this class is defined over the
interval $[-5, 5]$ and has the following form
\begin{equation}\label{fpinter}
f_j(x)=0.025(x-x_j^*)^2+sin^2((x-x_j^*)+(x-x_j^*)^2)+sin^2(x-x_j^*),
\end{equation}
where the global minimizer $x_j^*$, $j=1,...,100$, is chosen
randomly from interval $[-5, 5]$ and differently for the 100
functions of the class. Fig. \ref{fig.5} shows the graph of the
function no. 38 from the set of test functions (\ref{fpinter}) and
the trial points generated by the 12 methods while minimizing this
function, with accuracy $\varepsilon=10^{-4}(b-a)$. The global
minimum of the function, $f^*=0$, is attained at the point $x^*=
3.3611804993$. In Fig. \ref{fig.5}    the effects of the
acceleration techniques, the local tuning and the local
improvement, can be  clearly seen.

Table \ref{table5} shows the average numbers of trial points
generated by the six methods that do not use the derivative of the
objective functions, while Table \ref{table6} contains the results
of the six methods that use the derivative. In columns 2 and 4, the
values of the reliability parameter $r$ are given.   In Table
\ref{table5}, the asterisk denotes that in the algorithm $LT\_LI$
(for $\varepsilon=10^{-4}(b-a)$) the value r=1.3 has been used for
99 functions, and for the function no. 32   the value r=1.4 has been
applied.  Tables \ref{table5} and \ref{table6} confirm for the
second series of experiments  the same conclusions that have been
made with respect to the effects of the introduction of the
acceleration techniques for the first series of numerical tests.

\section{Conclusions}
\setcounter{equation}{0} In this paper,  there have been
considered two kinds of the one-dimensional global optimization
problems over a closed finite interval: (i) problems where the
objective function $f(x)$ satisfies the Lipschitz condition with a
constant $L$; (ii) problems where the first derivative of $f(x)$
satisfies the Lipschitz condition with a constant $M$.

Two general schemes describing   numerical methods for solving
both problems have been described. Six particular algorithms have
been presented for the case (i) and six algorithms for the case
(ii). In both cases, auxiliary functions constructed and
adaptively improved during the search have been used. In the case
(i), piece-wise linear functions have been described and
constructed. In the case (ii), smooth piece-wise quadratic
functions have been applied.

In the introduced methods, the Lipschitz constants $L$ and $M$ for
the objective function and its first derivative were either taken
as values known a priori or were dynamically estimated during the
search. A recent technique that adaptively estimates the local
Lipschitz constants over different zones of the search region was
used to accelerate the search in both cases (i) and (ii) together
with a newly introduced technique called the  local improvement.
Convergent conditions of the described twelve algorithms have been
studied.

The proposed local improvement technique is of a particular interest
due to the following  reasons. First, usually in the global
optimization methods the local search phases are separated from the
global ones. This means that it is necessary to introduce a rule
that  stops the global phase and starts the local one; then it stops
the local phase and starts the global one. It happens very often
that the global search and the local one are realized by different
algorithms and the global search is not able to use \textit{all}
evaluations of the objective function made during the local search
losing so an important information about the objective function that
has been already obtained.

The local improvement technique introduced in this paper does not
have this drawback  and allows the global search to use all the
information obtained during the local phases. In addition, it can
work both with and without  the derivatives and this is a valuable
asset when one solves the Lipschitz global optimization problems
because, clearly, Lipschitz functions can be non-differentiable.

Numerical experiments executed on 120 test problems taken from the
literature have shown quite a promising performance of the new
accelerating techniques.


\vspace{10mm}

\begin{thebibliography}{10}

\bibitem{Breiman&Cutler(1993)}
L.~Breiman and A.~Cutler.
\newblock A deterministic algorithm for global optimization.
\newblock {\em Mathematical Programming}, 58(1--3):179--199, 1993.

\bibitem{Calvin}
J.~M. Calvin.
\newblock An adaptive univariate global optimization algorithm and its
  convergence rate under the {W}iener measure.
\newblock {\em Informatica}, 22(4):471--488, 2011.

\bibitem{Calvin&Zilinskas}
J.~M. Calvin and A.~\u{Z}ilinskas.
\newblock One-dimensional global optimization for observations with noise.
\newblock {\em Computers and Mathematics with Applications}, 50(1--2):157--169,
  2005.

\bibitem{Casado:et:al.(2002)}
L.~G. Casado, I.~Garcia, and Ya.~D. Sergeyev.
\newblock Interval algorithms for finding the minimal root in a set of
  multiextremal non-differentiable one-dimensional functions.
\newblock {\em SIAM J. Scientific Computing}, 24(2):359--376, 2002.

\bibitem{Daponte:et:al.(1995)}
P.~Daponte, D.~Grimaldi, A.~Molinaro, and {Ya}.~D. Sergeyev.
\newblock An algorithm for finding the zero-crossing of time signals with
  {Lipschitzean} derivatives.
\newblock {\em Measurement}, 16(1):37--49, 1995.

\bibitem{Daponte:et:al.(1996)}
P.~Daponte, D.~Grimaldi, A.~Molinaro, and {Ya}.~D. Sergeyev.
\newblock Fast detection of the first zero-crossing in a measurement signal
  set.
\newblock {\em Measurement}, 19(1):29--39, 1996.

\bibitem{Floudas&Pardalos(1996)}
C.~A. Floudas and P.~M. Pardalos.
\newblock {\em State of the Art in Global Optimization}.
\newblock Kluwer Academic Publishers, Dordrecht, 1996.

\bibitem{Gergel(1992)}
V.~P. Gergel.
\newblock A global search algorithm using derivatives.
\newblock In Yu.I. Neymark, editor, {\em Systems Dynamics and Optimization},
  pages 161--178. N. Novgorod University Press, 1992.

\bibitem{Hamacher}
K.~Hamacher.
\newblock On stochastic global optimization of one-dimensional functions.
\newblock {\em Physica A: Statistical Mechanics and its Applications}, 354(15
  August 2005):547--557, 2005.

\bibitem{Hansen(1979)}
E.~R. Hansen.
\newblock Global optimization using interval analysis: The one-dimensional
  case.
\newblock {\em J. Optim. Theory Appl.}, 29(3):251--293, 1979.

\bibitem{Hansen&Jaumard(1992)}
P.~Hansen, B.~Jaumard, and H.~Lu.
\newblock Global optimization of univariate lipschitz functions: 1-2.
\newblock {\em Mathematical Programming}, 55:251--293, 1992.

\bibitem{Horst&Pardalos(1995)}
R.~Horst and P.~M. Pardalos, editors.
\newblock {\em Handbook of Global Optimization}, volume~1.
\newblock Kluwer Academic Publishers, Dordrecht, 1995.

\bibitem{Horst&Tuy(1996)}
R.~Horst and H.~Tuy.
\newblock {\em Global Optimization -- Deterministic Approaches}.
\newblock Springer--Verlag, Berlin, 1996.

\bibitem{kn:theory}
D.~E. Johnson.
\newblock {\em Introduction to Filter Theory}.
\newblock Prentice Hall Inc., New Jersey, 1976.

\bibitem{kn:KABA89}
D.~Kalra and A.~H. Barr.
\newblock Guaranteed ray intersections with implicit surface.
\newblock {\em Computer Graphics}, 23(3):297--306, 1989.

\bibitem{Kvasov&Sergeyev(2009)}
D.~E. Kvasov and {Ya.}~D. Sergeyev.
\newblock A univariate global search working with a set of {Lipschitz}
  constants for the first derivative.
\newblock {\em Optimization Letters}, 3(2):303--318, 2009.

\bibitem{kn:filters}
H.~Y.-F. Lam.
\newblock {\em Analog and Digital Filters-Design and Realization}.
\newblock Prentice Hall Inc., New Jersey, 1979.

\bibitem{Locatelli}
M.~Locatelli.
\newblock Bayesian algorithms for one-dimensional global optimization.
\newblock {\em J. Global Optimization}, 10(1):57--76, 2001.

\bibitem{Mladineo(1992)}
R.~Mladineo.
\newblock Convergence rates of a global optimization algorithm.
\newblock {\em Mathematical Programming}, 54:223--232, 1992.

\bibitem{Molinaro&Sergeyev(2001b)}
A.~Molinaro and {Ya}.~D. Sergeyev.
\newblock Finding the minimal root of an equation with the multiextremal and
  nondifferentiable left-hand part.
\newblock {\em Numerical Algorithms}, 28(1--4):255--272, 2001.

\bibitem{Pinter(1996)}
J.~D. {Pint\'{e}r}.
\newblock {\em Global Optimization in Action (Continuous and Lipschitz
  Optimization: Algorithms, Implementations and Applications)}.
\newblock Kluwer Academic Publishers, Dordrecht, 1996.

\bibitem{Pinter(2002)}
J.~D. {Pint\'{e}r}.
\newblock Global optimization: software, test problems, and applications.
\newblock In P.~M. Pardalos and H.~E. Romeijn, editors, {\em Handbook of Global
  Optimization}, volume~2, pages 515--569. Kluwer Academic Publishers,
  Dordrecht, 2002.

\bibitem{Piya(1972)}
S.~A. Piyavskii.
\newblock An algorithm for finding the absolute extremum of a function.
\newblock {\em Com. Maths. Math. Phys.}, 12:57--67, 1972.

\bibitem{Sergeyev(1995b)}
{Ya.}~D. Sergeyev.
\newblock A one-dimensional deterministic global minimization algorithm.
\newblock {\em Comput. Math. Math. Phys.}, 35(5):705--717, 1995.

\bibitem{Sergeyev(1998)}
{Ya.}~D. Sergeyev.
\newblock Global one-dimensional optimization using smooth auxiliary functions.
\newblock {\em Mathematical Programming}, 81(1):127--146, 1998.

\bibitem{Sergeyev(2006)}
{Ya.}~D. Sergeyev.
\newblock Univariate global optimization with multiextremal non-differentiable
  constraints without penalty functions.
\newblock {\em Computational Optimization and Applications}, 34(2):229--248,
  2006.

\bibitem{Sergeyev:et:al.(1999)}
{Ya.}~D. Sergeyev, P.~Daponte, D.~Grimaldi, and A.~Molinaro.
\newblock Two methods for solving optimization problems arising in electronic
  measurements and electrical engineering.
\newblock {\em SIAM J. Optim.}, 10(1):1--21, 1999.

\bibitem{Sergeyev:et:al(2001)}
{Ya.}~D. Sergeyev, D.~Famularo, and P.~Pugliese.
\newblock Index branch-and-bound algorithm for lipschitz univariate global
  optimization with multiextremal constraints.
\newblock {\em J. Global Optimization}, 21(3):317--341, 2001.

\bibitem{Sergeyev&Grishagin(1994)}
{Ya.}~D. Sergeyev and V.A. Grishagin.
\newblock A parallel algorithm for finding the global minimum of univariate
  functions.
\newblock {\em J. Optimization Theory and Applications}, 80(3):513--536, 1994.

\bibitem{Sergeyev&Kvasov(2008)}
{Ya.}~D. Sergeyev and D.~E. Kvasov.
\newblock {\em Diagonal Global Optimization Methods}.
\newblock {FizMatLit}, Moscow, 2008.
\newblock In Russian.

\bibitem{Sergeyev:et:al.(2007)}
{Ya.}~D. Sergeyev, D.~E. Kvasov, and F.~M.~H. Khalaf.
\newblock A one-dimensional local tuning algorithm for solving {GO} problems
  with partially defined constraints.
\newblock {\em Optimization Letters}, 1(1):85--99, 1995.

\bibitem{Sergeyev&Markin(1995)}
{Ya.}~D. Sergeyev and D.~L. Markin.
\newblock An algorithm for solving global optimization problems with nonlinear
  constraints.
\newblock {\em J. Global Optim.}, 7(4):407--419, 1995.

\bibitem{Strongin&Sergeyev(2000)}
R.~G. Strongin and {Ya.}~D. Sergeyev.
\newblock {\em Global optimization with non-convex constraints: {Sequential}
  and parallel algorithms}.
\newblock Kluwer Academic Publishers, Dordrecht, 2000.

\bibitem{Torn&Zilinskas(1989)}
A.~{T\"orn} and A.~{\v{Z}ilinskas}.
\newblock {\em Global Optimization}, volume 350 of {\em Lecture Notes in
  Computer Science}.
\newblock Springer--Verlag, Berlin, 1989.

\bibitem{Zilinskas}
A.~\u{Z}ilinskas.
\newblock Optimization of one-dimensional multimodal functions: Algorithm {AS}
  133.
\newblock {\em Applied Statistics}, 23:367--375, 1978.

\bibitem{Wolfe}
M.~A. Wolfe.
\newblock On first zero crossing points.
\newblock {\em Applied Mathematics and Computation}, 150:467--479, 2004.

\end{thebibliography}
\end{document}